\documentclass[12pt]{amsart}
\usepackage{amssymb,amsmath,a4wide}%,showtags}

\newtheorem{theorem}{Theorem}[section]
\newtheorem{lemma}[theorem]{Lemma}

\newtheorem{remark}{\it Remark\/}

\newcommand{\pa}{\partial}
\newcommand{\ep}{\varepsilon}
\newcommand{\RR}{{\mathbb{R}}}
\renewcommand{\div}{\mathop{{\rm div}}}
\newcommand{\rest}[1]{\rule[-8pt]{.5pt}{14pt}_{\;#1}}

\begin{document}

\title[Loewner-Nirenberg problem]{Boundary behavior in the Loewner-Nirenberg problem}%
\author{Satyanad Kichenassamy}%
\address{Laboratoire de Math\'ematiques (UMR 6056), CNRS \&\ Universit\'e
de Reims Champagne-Ardenne, Moulin de la Housse, B.P. 1039,
F-51687
Reims Cedex 2\\France}%
\email{satyanad.kichenassamy@univ-reims.fr}%

\thanks{}%

\date{May 10, 2004}%
%\dedicatory{}%
%\commby{}%
% ----------------------------------------------------------------
\vskip 1em

\hfill{\bf Appeared in:} \emph{Journal of Functional Analysis},
{\bf 222} : 1 (2005) 98--113.

\vskip 1em
\begin{abstract}
Let $\Omega\subset\RR^n$ be a bounded domain of class
$C^{2+\alpha}$, $0<\alpha<1$. We show that if $n\geq 3$ and
$u_\Omega$ is the maximal solution of equation $\Delta u =
n(n-2)u^{(n+2)/(n-2)}$ in $\Omega$, then the hyperbolic radius
$v_\Omega=u_\Omega^{-2/(n-2)}$ is of class $C^{2+\alpha}$ up to
the boundary. The argument rests on a reduction to a nonlinear
Fuchsian elliptic PDE.
\end{abstract}
\maketitle
%\tableofcontents
% ----------------------------------------------------------------

\section{Introduction}
\label{sec:intro}

\subsection{Main result}
Let $\Omega\subset\RR^n$, $n\geq 3$, be a bounded domain of class
$C^{2+\alpha}$, where $0<\alpha<1$. Consider the Loewner-Nirenberg
equation in the form
\begin{equation}\label{eq:LN}
-\Delta u + n(n-2)u^{\frac{n+2}{n-2}}=0.
\end{equation}
It is known \cite{LN} that this equation admits a maximal solution
$u_\Omega$, which is positive and smooth inside $\Omega$; it is
the limit of the increasing sequence $(u_m)_{m\geq 1}$ of
solutions of (\ref{eq:LN}) which are equal to $m$ on the boundary.
The \emph{hyperbolic radius} of $\Omega$ is the function
\[
v_\Omega:=u_\Omega^{-2/(n-2)};
\]
it vanishes on $\pa\Omega$. Let $d(x)$ denote the distance of $x$
to $\pa\Omega$. It is of class $C^{2+\alpha}$ near $\pa\Omega$. We
prove
\begin{theorem}\label{th:1}
If $\Omega$ is of class $C^{2+\alpha}$, then $v_\Omega\in
C^{2+\alpha}(\overline\Omega)$, and
\[ v_\Omega(x) = 2d(x)-d(x)^2[H(x)+o(1)]
\]
as $d(x)\to 0$, where $H(x)$ is the mean curvature at the point of
$\pa\Omega$ closest to $x$.
\end{theorem}
This result is optimal, since $H$ is of class $C^\alpha$ on the
boundary. It follows from Theorem \ref{th:1} that $v_\Omega$ is a
\emph{classical solution} of
\[v_\Omega\Delta v_\Omega=\frac{n}2(|\nabla v_\Omega|^2-4),
\]
even though $u_\Omega$ cannot be interpreted as a weak solution of
(\ref{eq:LN}), insofar as $u^{\frac{n+2}{n-2}}\sim
(2d)^{-1-n/2}\not\in L^1(\Omega)$.

\subsection{Motivation}

The main reasons for studying $u_\Omega$ are as follows.
\begin{itemize}
\item $u_\Omega$ dominates all classical solutions, and therefore
provides a uniform interior bound, independent of boundary data
(see \cite{K,O,LN}).
\item The metric
\[ v_\Omega^{-2}(dx_1^2+\cdots+dx_n^2)
\]
is complete, and has constant negative scalar curvature; it
therefore generalizes the Poincar\'e metric on the unit disc and
provides an intrinsic geometry on $\Omega$. Furthermore, equation
(\ref{eq:LN}) admits a partial conformal invariance property. This
was the motivation of Loewner and Nirenberg \cite{LN}.
\item The minima of $v_\Omega$, known as hyperbolic centers, are
close to the points of concentration arising in several
variational problems of recent interest, see \cite{BF}.
\end{itemize}
The numerical computation of $v_\Omega$ proceeds by computing the
solution of the Dirichlet problem for (\ref{eq:LN}) on a set of
the form $\{d(x)>h\}$, where $h$ is small, and the Dirichlet data
are given by the boundary asymptotics of $u_\Omega$.

Now, earlier results on the boundary behavior of $u_\Omega$,
summarized below, yield
\begin{equation}\label{eq:asy-v}
v_\Omega=2d+O(d^2)\text{ and }|\nabla v_\Omega|\to 2
\end{equation}
as $d(x)\to 0$. Motivated by this, Bandle and Flucher conjectured
Theorem \ref{th:1} (\cite[p.~204]{BF}).

In two dimensions, the hyperbolic radius is defined by
$v_\Omega=\exp(-u_\Omega)$, where $u_\Omega$ solves the Liouville
equation
\[
-\Delta u + 4e^{2u}=0.
\]
For background information on the two-dimensional case, see
\cite{BF,B}; one should add that, if $\Omega$ is simply connected,
the hyperbolic radius coincides with the conformal or mapping
radius, and with the harmonic radius. The analogue of Theorem
\ref{th:1} for the Liouville equation is proved in
\cite{SK-note,SK-be}.

\subsection{Earlier results}

Loewner and Nirenberg showed the existence of $u_\Omega$ and
proved that
\begin{itemize}
  \item if
   $\Omega\subset\Omega'$, then any classical solution in $\Omega'$
   restricts to a classical solution in $\Omega$, so that
   \begin{equation}\label{eq:monot}
   u_{\Omega'}\leq u_\Omega;
  \end{equation}
  \item $u\sim (2d)^{1-\frac n2}$ as $d\to 0$.
\end{itemize}

It follows from \cite{KN,LM} that $u=(2d)^{1-n/2}(1+O(d))$ as
$d\to 0$.

It follows from \cite[pp.~95--96]{BE} and \cite{BM} that
\[ |\nabla u_\Omega|(2d)^{n/2}\to (n-2).
\]
From this information, equation (\ref{eq:asy-v}) follows.

There is an extensive literature on the issue of boundary blow-up,
see  \cite{BE,BF,BM,BMd,BP,K,KN,LM,LN,MV,O} and their references for
details.

\subsection{Method of proof}

We begin by performing a Fuchsian reduction, that is, we introduce
the degenerate equation solved by a renormalized unknown, which
governs the higher-order asymptotics of the solution; in this
case, a convenient renormalized unknown is
\[
w:=(v_\Omega-2d)/d^2.
\]
It follows from general arguments, see the overview in \cite{SK},
that the equation for $w$ has a very special structure: the
coefficient of the derivatives of order $k$ is divisible by $d^k$
for $k=0$, 1 and 2, and the nonlinear terms all contain a factor
of $d$. Such an equation is said to be \emph{Fuchsian}; the
regularity properties of solutions of such equations are discussed
in section \ref{sec:f}.

In the present case, one finds
\begin{equation}\label{eq:fr}
\frac {2v^{n/2}}{n-2}\{ -\Delta
u_\Omega+n(n-2)u_\Omega^{(n+2)/(n-2)}\} =Lw+2\Delta d-M_w(w),
\end{equation}
where
\[
L:=d^2\Delta+(4-n)d\nabla d\cdot\nabla+(2-2n),
\]
and $M_w$ is a linear operator with $w$-dependent coefficients,
defined by
\[
M_w(f):=\frac{nd^2}{2(2+dw)} [2f\nabla d\cdot\nabla w+d\nabla
w\cdot\nabla f]-2df\Delta d.
\]
The proof now consists in a careful bootstrap argument in which
better and better information on $w$ results in better and better
properties of the degenerate linear operator $L-M_w$. A key step
is the inversion of the analogue of $L$ in the half-space, which
plays the role of the Laplacian in the usual Schauder theory.

Equation (\ref{eq:fr}) needs only to be studied in the
neighborhood of the boundary. Let us therefore introduce
$C^{2+\alpha}$ thin domains $\Omega_\delta=\{0<d<\delta\}$, such
that $d\in C^{2+\alpha}(\overline\Omega_\delta)$, and
$\pa\Omega_\delta=\pa\Omega\cup\Gamma$ consists of two
hypersurfaces of class $C^{2+\alpha}$.

We will use the spaces
$C_\sharp^{k+\alpha}(\overline\Omega_\delta)$, consisting of
functions $w$ such that, for every $0\leq j\leq k$, $d^jw\in
C^{j+\alpha}(\overline\Omega_\delta)$. We write
\[ \|u\|_{C_\sharp^{k+\alpha}(\overline\Omega_\delta)}
:=\sum_{j=0}^k \|d^ju\|_{C^{j+\alpha}(\overline\Omega_\delta)}.
\]

The proof proceeds in five steps, corresponding to the following
theorems:
\begin{theorem}\label{th:I}
$w$ and $d^2\nabla w$ are bounded near $\pa\Omega$.
\end{theorem}
\begin{theorem}\label{th:II}
$d\nabla w$ and $M_w(w)/d$ are bounded near $\pa\Omega$.
\end{theorem}
\begin{theorem}\label{th:III}
If $\delta$ is sufficiently small, there is a $w_0\in
C_\sharp^{2+\alpha}(\overline\Omega_\delta)$ such that
\begin{equation}\label{eq:Lw0}
Lw_0+2\Delta d=0
\end{equation}
in $\Omega_\delta$. Furthermore
\begin{equation}\label{eq:H}
w_0\,\rest{\pa\Omega}=-H,
\end{equation}
where $H=-(\Delta d)/(n-1)$ is the mean curvature of the boundary.
\end{theorem}
\begin{theorem}\label{th:IV}
Near the boundary,
\[
\tilde{w}:=w-w_0=O(d).
\]
\end{theorem}
\begin{theorem}\label{th:V}
$\tilde{w}$ belongs to
$C_\sharp^{2+\alpha}(\overline\Omega_\delta)$.
\end{theorem}
Since $w=w_0+\tilde{w}$, we obtain $w\in
C_\sharp^{2+\alpha}(\overline\Omega_\delta)$, and since
$\tilde{w}=O(d)$, $w\,\rest{\pa\Omega}$ is equal to $-H$. This
completes the proof.

\subsection{Organization of the paper}

Section \ref{sec:f} recalls some results on linear Fuchsian PDE
from \cite{SK-be}.

Section \ref{sec:comp1} gives the proofs of Theorems \ref{th:I}
and  \ref{th:II}.

Section \ref{sec:prel} introduces a convenient coordinate system
near the boundary, and solves $Lw_1+2\Delta d=O(d^\alpha)$, as a
preparation for the proof of Theorem \ref{th:III}. It rests on the
analysis of the analogue of $L$ for the half-space.

Section \ref{sec:pfIII} proves Theorem \ref{th:III}.

Section \ref{sec:comp2} proves Theorem \ref{th:IV}.

Section \ref{sec:synth} gives the proof of Theorem \ref{th:V},
thus completing the proof of the main result.

% ----------------------------------------------------------------

\section{Background results on Fuchsian PDE}
\label{sec:f}

Let $(a^{ij})\in C^\alpha(\overline\Omega_\delta)$ be uniformly
elliptic. Recall that $\delta$ is chosen small enough so that
$d\in C^{2+\alpha}(\overline\Omega_\delta)$.

An operator $A$ is said to be of type (I) if it has the form
\[
A:=\pa_i(d^2a^{ij}(x)\pa_{j})+db^i(x)\pa_i+c(x),
\]
where $b^i$ and $c$ belong to $L^\infty(\Omega_\delta)$.

It is said to be of type (II) if it has the form
\[
A:=d^2a^{ij}(x)\pa_{ij}+db^i(x)\pa_i+c(x),
\]
where $b^i$ and $c$ belong to $C^\alpha(\overline\Omega_\delta)$.

Operator $L$, which may be written $\div(d^2\nabla)+(2-n)d\nabla
d\cdot\nabla+(2-2n)$, is of type (I) as well as (II).

The results we will need are the following.
\begin{theorem}\label{th:FIa}
If
\begin{enumerate}
  \item $A$ is of type (I), and
  \item $Af$ and $f$ are in $L^\infty(\Omega_\delta)$,
\end{enumerate}
then $df$ and $d^2\nabla f$ belong to
$C^\alpha(\overline\Omega_{\delta'})$ for $\delta'<\delta$, and
$d\nabla f$ is bounded near $\pa\Omega$.
\end{theorem}
This is proved in two dimensions in Theorem 5.1 of \cite{SK-be};
the proof applies without modification in $n$ dimensions.
\begin{remark}
If we also know that $f$ is, say of class $C^{2+\alpha}$ on
$\Gamma$, the conclusion holds on all of $\overline\Omega_\delta$;
a similar remark applies to the next two theorems as well.
\end{remark}
\begin{theorem}\label{th:FIb}
If
\begin{enumerate}
  \item $A$ is of type (I), and
  \item $Af$ and $f$ are $O(d^\alpha)$ as $d\to 0$,
\end{enumerate}
then $f\in C_\sharp^{1+\alpha}(\overline\Omega_{\delta'})$ for
$\delta'<\delta$.
\end{theorem}
This corresponds to Theorem 5.2 in \cite{SK-be}: in the latter
paper, it is assumed that $Af=O(d)$, and that $n=2$, but the proof
proceeds \emph{verbatim} for any $n$, if one only knows that
$Af=O(d^\alpha)$.
\begin{theorem}\label{th:FIIa}
If
\begin{enumerate}
  \item $A$ is of type (II),
  \item $Af\in C^\alpha(\overline\Omega_\delta)$,
  \item $f\in C_\sharp^{1+\alpha}(\overline\Omega_\delta)$,
\end{enumerate}
then $f\in C_\sharp^{2+\alpha}(\overline\Omega_{\delta'})$ for
$\delta'<\delta$.
\end{theorem}
\begin{proof}
The assumptions ensure that $a^{ij}\pa_{ij}(d^2f)$ is
H\"older-continuous and that $f$ is bounded; $d^2f$ therefore
solves a Dirichlet problem to which the Schauder estimates apply
near $\pa\Omega$. Therefore $d^2f$ is of class $C^{2+\alpha}$ up
to the boundary. Since we already know that $f\in
C^\alpha(\overline\Omega_\delta)$ and $df$ is of class
$C^{1+\alpha}(\overline\Omega_\delta)$, we have indeed $f$ of
class $C_\sharp^{2+\alpha}(\overline\Omega_{\delta'})$ for
$\delta'<\delta$.
\end{proof}

% ----------------------------------------------------------------

\section{First comparison argument and proof of Theorems \ref{th:I} and
\ref{th:II}} \label{sec:comp1}

\subsection{Proof of Theorem \ref{th:I}}

We give a self-contained proof for the convenience of the reader.
It could be slightly shortened if one starts from the information
$v_\Omega\sim 2d$ from \cite{LN}.

Since $\pa\Omega$ is $C^{2+\alpha}$, it satisfies a uniform
interior and exterior sphere condition, and there is a positive
$r_0$ such that any $P\in\Omega$ such that $d(P)\leq r_0$ admits a
unique nearest point $Q$ on the boundary, and such that there are
two points $C$ and $C'$ on the line determined by $P$ and $Q$,
such that
\[B_{r_0}(C)\subset\Omega\subset\RR^n\setminus B_{r_0}(C'),\]
these two balls being tangent to $\pa\Omega$ at $Q$. We now define
two functions $u_i$ and $u_e$. Let
\[u_i(M)=(r_0-\frac{CM^2}{r_0})^{1-n/2}\text{ and }
u_e(M)=(\frac{C'M^2}{r_0}-r_0)^{1-n/2}.
\]
$u_i$ and $u_e$ are solutions of equation (\ref{eq:LN}) in
$B_{r_0}(C)$ and $\RR\setminus B_{r_0}(C')$ respectively.

If we replace $r_0$ by $r_0-\varepsilon$ in the definition of
$u_e$, we obtain a classical solution of (\ref{eq:LN}) in
$\Omega$, which is therefore dominated by $u_\Omega$. It follows
that
\[ u_e\leq u_\Omega\text{ in }\Omega.\]
The monotonicity property (\ref{eq:monot}) yields
\[ u_\Omega\leq u_i\text{ in }B_{r_0}(C).\]

In particular, the inequality
\[ u_e(M)\leq u_\Omega(M)\leq u_i(M)\]
holds if $M$ lies on the semi-open segment $[P,Q)$. Since $Q$ is
then also the point of the boundary closest to $M$, we have
$QM=d(M)$, $CM=r_0-d$ and $C'M=r_0+d$; it follows that
\[(2d+\frac{d^2}{r_0})^{1-n/2}\leq u_\Omega(M)
\leq (2d-\frac{d^2}{r_0})^{1-n/2}.
\]
Since $u_\Omega=(2d+d^2w)^{1-n/2}$, it follows that
\[ |w|\leq \frac1{r_0}\text{ if }d\leq r_0.  \]

Next, consider $P\in\Omega$ such that $d(P)=2\sigma$, with
$3\sigma<r_0$. For $x$ in the closed unit ball $\overline B_1$,
let
\[  P_\sigma :=P+\sigma x;\quad
    u_\sigma(x) := \sigma^{(n-2)/2}u(P_\sigma).
\]
One checks that $u_\sigma$ is a classical solution of
(\ref{eq:LN}) in $\overline B_1$. Since $d\mapsto 2d\pm \frac
1{r_0}d^2$ is increasing for $d<r_0$, and $d(P_\sigma)$ varies
between $\sigma$ and $3\sigma$ if $x$ varies in $\overline B_1$,
we have
\[(6+\frac{9\sigma}{r_0})^{1-n/2}\leq u_\sigma(M)
\leq (2-\frac{\sigma}{r_0})^{1-n/2}.
\]
This provides a uniform bound for $u_\sigma$ on $B_1$. Applying
interior regularity estimates as in \cite{SK-th,BM}, we find that
$\nabla u_\sigma$ is uniformly bounded for $x=0$. Recalling that
$\sigma=\frac12 d(P)$, we find that
\[ d^{\frac n2-1}u\text{ and }d^{\frac n2}\nabla u\text{ are
bounded near }\pa\Omega.
\]
It follows that $u^{-n/(n-2)}=O(d^{n/2})$, and since
$d^2w=-2d+u^{-2/(n-2)}$, we have
\[ d^2\nabla w=-2(1+dw)\nabla d-\frac 2{n-2}u^{-n/(n-2)}\nabla u,
\]
hence $d^2\nabla w$ is bounded near $\pa\Omega$. This completes
the proof of Theorem \ref{th:I}.

\subsection{Proof of Theorem \ref{th:II}}

Since $w$ and $d^2\nabla w$ are bounded near the boundary, it
follows that $L-M_w$ is an operator of type (I) near $\pa\Omega$.
Theorem \ref{th:FIa} now ensures that $d\nabla w$ is bounded near
$\pa\Omega$. Going back to the definition of $M_w$, we find that
$M_w(w)=O(d)$ near the boundary.

This completes the proof of Theorem \ref{th:II}.

Note that Theorem \ref{th:FIa} gives in addition that $dw$ and
$d^2\nabla w$ are of class $C^\alpha$ near the boundary.

At this stage, we have proved that
\[ Lw+2\Delta d=O(d).\]

% ----------------------------------------------------------------

\section{The $(Y,T)$ coordinates and the model operator $L_0$}
\label{sec:prel}

We write henceforth $u$ and $v$ for $u_\Omega$ and $v_\Omega$
respectively.

\subsection{Local coordinates near a point of the boundary}

Since $\pa\Omega$ is compact, there is a positive $r_0$ such that
in any ball of radius $r_0$ centered at a point of $\pa\Omega$,
one may introduce a coordinate system $(Y,T)$ in which $T=d$ is
the last coordinate. It will be convenient to assume that the
domain of this coordinate system contains a set of the form
\[ 0<T<\theta\text{ and }|Y_j|<\theta\text{ for }j\leq n-1.
\]
Let $\pa_j=\pa_{x_j}$, and write $d_n$ and $d_j$ for $\pa d/\pa
x_n$ and $\pa d/\pa x_j$ respectively. Primes denote derivatives
with respect to the $Y$ variables: $\pa'_j=\pa_{Y_j}$,
$\nabla'=\nabla_Y$, $\Delta'=\sum_{j<n}\pa_j^{\prime 2}$, etc.;
note that this definition of $\Delta'$ differs from the one in
\cite{SK-be}.  We write $\tilde\nabla d=(d_1,\dots,d_{n-1})$.
Recall that $|\nabla d|=1$. We let throughout
\[
D=T\pa_T.
\]
The transformation formulae are
\begin{eqnarray*}
    T & =&  d(x_1,\dots,x_n); \qquad Y_j = x_j\text{ for }j<n;\\
    \pa_n &=& d_n\pa_T; \qquad \pa_j=d_j\pa_T+\pa_j'.
\end{eqnarray*}
We recall that if the coordinate axes are such that the origin is
on $\pa\Omega$, and the $x_n$ axis points in the direction of the
inward normal, then the origin $(Y=0, T=0)$ is the point of
$\pa\Omega$ closest to $(0,T)$ for $T$ small; the line $T\mapsto
(0,T)$ is the normal to the boundary; furthermore, if
$\kappa_i(Y)$ are the principal curvatures of the boundary, we
have, for $Y=0$,
\[ \quad \Delta d(0,T)
     =-\sum_1^{n-1}\frac{\kappa_i(0)}{1-T\kappa_i(0)},
\]
see \cite[\S 14.6]{GT}. It follows that, on the boundary,
\[ 2\Delta d=(2-2n)H,
\]
where $H$ is the mean curvature of $\pa\Omega$.

We further have
\begin{eqnarray*}
    d\nabla d\cdot\nabla w &=& (D+T\tilde\nabla d\cdot\nabla')w \\
    |\nabla w|^2 &=& w_T^2+|\nabla'
    w|^2+2w_T\tilde\nabla d\cdot\nabla'w\\
    \Delta w &=&
    w_{TT}+\Delta'w+2\tilde\nabla d\cdot\nabla'w_T+w_T\Delta d.
\end{eqnarray*}
It follows that
\[ Lw=L_0w+L_1w,\]
where
\[ L_0w=(D+2)(D+1-n)w+T^2\Delta'w,
\]
and
\[ L_1w=(4-n)\tilde\nabla d\cdot\nabla'(Tw)
        +2T\tilde\nabla d\cdot\nabla'(Dw)+T(Dw)\Delta d.
\]

\subsection{Solution of $Lf=k+O(d^\alpha)$}

Let $C^\alpha_{\text{per}}$ denote the space of functions
$k(Y,T)\in C^\alpha(0\leq T\leq \theta)$ which satisfy
$k(Y_j+2\theta,T)=k(Y_j,T)$ for $1\leq j\leq n-1$ . We prove the
following theorem.
\begin{theorem}\label{th:L0}
Let $\theta>0$, and $k(Y,T)$ of class $C^\alpha_{\text{per}}$ Then
there is a function $f$ such that
\begin{enumerate}
  \item $L_0f=k+O(d^\alpha)$,
  \item $f$ is of class $C_\sharp^{2+\alpha}(0\leq T\leq\theta)$,
  \item $f(Y,0)=k(Y,0)/(2-2n)$ and
  \item $L_1f=O(d^\alpha)$.
\end{enumerate}
\end{theorem}
\begin{proof}
Let
\[ L'_0=(D+2)(D-1)+T^2\Delta'=L_0+(n-2)(D+2).
\]
We first solve the equation $L'_0f_0=k$, as in \cite{SK-be}.
\begin{lemma} There is a bounded linear operator $G$ from
$C^\alpha_{\text{per}}$ to $C_\sharp^{2+\alpha}(0\leq
T\leq\theta)$ such that $f_0:=G[k]$ verifies
\begin{enumerate}
  \item $L'_0f_0=k$,
  \item $f_0$ is of class $C_\sharp^{2+\alpha}(0\leq T\leq\theta)$,
  \item $f_0(Y,0)+k(Y,0)/2=0$, $Df_0(Y,0)=0$ and
  \item $L_1f_0=O(d^\alpha)$.
\end{enumerate}
\end{lemma}
\begin{proof}
One first constructs $\tilde k$ such that $(D-1)\tilde k=-k$, and
$\tilde k$ and $D\tilde k$ are both $C^\alpha$ up to $T=0$. One
may take
\[ \tilde{k}=\int_1^\infty
F_1[k](T\sigma, Y)\frac{d\sigma}{\sigma^2}.
\]
where $F_1$ is an extension operator, so that $F_1[k]=k$ for
$T\leq \theta$.

One checks that $\tilde k=k$ for $T=0$.

One then solves $(\pa_{TT}+\Delta')h+\tilde k=0$ with periodic
boundary conditions, of period $2\theta$, in each of the $Y_j$,
and $h(Y,0)=h_T(Y,\theta)=0$; this yields
\[ h \text{ is of class }C^{2+\alpha}(0\leq T\leq \theta)
\]
by the Schauder estimates. In particular, $h_T$ is continuous up
to $T=0$, and $Dh=0$ for $T=0$ and $T=\theta$.

Since $h=0$ for $T=0$, we also have $\Delta' h=0$ for $T=0$. The
equation for $h$ therefore gives
\[h_{TT}=-\tilde k=-k\text{ for }T=0.
\]

In addition,
\[(\pa_{TT}+\Delta')Dh=D(\pa_{TT}+\Delta')h+2h_{TT}=k-\tilde
k+2h_{TT},
\]
which is $C^\alpha$. Since, on the other hand, $Dh$ is of class
$C^1$ and $Dh=0$ for $T=0$ and $T=\theta$, we conclude, using
again the Schauder estimates, that
\[ Dh \text{ is of class }C^{2+\alpha}(0\leq T\leq \theta).
\]

We now define $f_0$ by
\begin{equation}\label{eq:G}
f_0:=T^{-2}(D-1)h
        =\pa_T\left(\frac h T\right)
        =\int_0^1 \sigma h_{TT}(Y,T\sigma)\,d\sigma.
\end{equation}
Since $f_0$ is itself uniquely determined by $h$, itself defined
in terms of $k$ we define a map $G$ by
\[f_0=G[k].
\]
A direct computation yields $L'_0f_0=k$, see \cite[section
6.1]{SK-be}.

Let us now consider the regularity of $f_0$ up to $\pa\Omega$, and
the values of $f_0$ and its derivatives on $\pa\Omega$.

Consider $g_0:=T^2f_0$. Since $g_0=(D-1)h\in C^{2+\alpha}(0\leq
T\leq \theta)$ and vanishes for $T=0$, we have $g_0=\int_0^1
g_{0T}(Y,T\sigma)Td\sigma$. It follows that
\[Tf_0(Y,T)=\int_0^1 g_{0T}(Y,T\sigma)d\sigma
\in C^{1+\alpha}(0\leq T\leq \theta).
\]
Since, on the other hand, $G[k]=\int_0^1 \sigma
h_{TT}(Y,T\sigma)\,d\sigma$, we find $f_0\in C^\alpha(0\leq T\leq
\theta)$, and
\[f_0(Y,0)=\frac12 h_{TT}(Y,0)=-\frac12 k(Y,0).
\]
We therefore have
\[f_0 \text{ is of class }C_\sharp^{2+\alpha}(0\leq T\leq \theta).
\]
Since
\[ (D+2)f_0 =T^{-2}D(D-1)h= h_{TT},
\]
we find $Df_0(Y,0)=h_{TT}(Y,0)-2f_0(Y,0)=0$. By differentiation
with respect to the $Y$ variables, we obtain that $\tilde\nabla
d\cdot\nabla'(Tf_0)$ is of class $C^\alpha$ and vanishes for
$T=0$. The same is true of $T(Df_0)\Delta d$. Similarly,
\[
2T\tilde\nabla d\cdot\nabla'Df_0 = 2\tilde\nabla
d\cdot\nabla'[\pa_T(T^2f_0)-2Tf_0]
\]
is of class $C^\alpha$, and vanishes for $T=0$ because this is
already the case for $TDf_0$. It follows that $L_1f_0$ is a
$C^\alpha$ function which vanishes for $T=0$; it is therefore
$O(d^\alpha)$ as desired.
\end{proof}

We are now ready to prove Theorem \ref{th:L0}. Let $a$ be a
constant, and $f=G[ak]$. We therefore have $L'_0f=ak$, and, for
$T=0$, $f=-\frac 12 ak$. Since $L_1f\in C^\alpha$, and $L_1f$ and
$Df$ both vanish for $T=0$, it follows that, for $T=0$,
\[
Lf-k=(L'_0-(n-2)(D+2)+L_1)f-k=[a+(n-2)a-1]k.
\]
Taking $a=1/(n-1)$, we find that $f$ has the announced properties.
\end{proof}

% ----------------------------------------------------------------

\section{Construction of $w_0$ and proof of Theorem \ref{th:III}}
\label{sec:pfIII}

\subsection{Solution of $Lw_0=g$}
\label{sec:consL}

Let us now consider a function $g$ of class
$C^\alpha(\overline\Omega_\delta)$.

Recall that there is a positive $r_0<\delta$ such that any ball of
radius $r_0$, centered at a point of the boundary, is contained in
a domain in which we have a system of coordinates of the type
$(Y,T)$. Let us cover (a neighborhood of) $\pa\Omega$ by a finite
number of balls $(V_\lambda)_{\lambda\in\Lambda}$ of radius
$r_1<r_0$ and centers on $\pa\Omega$, and consider the balls
$(U_\lambda)_{\lambda\in\Lambda}$ of radius $r_0$ with the same
centers. Thus, we may assume that every $U_\lambda$ is associated
with a coordinate system $(Y_\lambda,T_\lambda)$ of the type
considered in section \ref{sec:prel}; taking $r_1$ smaller if
necessary, we may also assume that $\overline{V}_\lambda\subset
Q_\lambda\subset U_\lambda$, where $Q_\lambda$ has the form
\[Q_\lambda:=\{(Y_{\lambda,1,\dots,}Y_{\lambda,n-1},T_\lambda) : 0\leq Y_{\lambda,
j}\leq\theta\text{ for every $j$, and }0<T_\lambda<\theta\}.
\]
Consider a smooth partition of unity $(\varphi_\lambda)$ and
smooth functions $(\Phi_\lambda)$, such that
\begin{enumerate}
  \item $\sum_{\lambda\in\Lambda}\varphi_\lambda=1$ near $\pa\Omega$;
  \item supp\,$\varphi_\lambda\subset V_\lambda$;
  \item supp\,$\Phi_\lambda\subset U_\lambda\cap\{T<\theta\}$;
  \item $\Phi_\lambda=1$ on $V_\lambda$.
\end{enumerate}
In particular, $\Phi_\lambda\varphi_\lambda=\varphi_\lambda$.

The function $g\varphi_\lambda$ is of class $C^\alpha(\overline
Q_\lambda)$; it may be extended by successive reflections to an
element of $C^\alpha_{\text{per}}$, with period $2\theta$ in the
$Y_\lambda$ variables; this extension will be denoted by the same
symbol for simplicity.

Let us apply Theorem \ref{th:L0}, and consider, for every
$\lambda$, the function $w_\lambda:=G[g\varphi_\lambda/(n-1)]$. We
have
\[
Lw_\lambda=g\varphi_\lambda+R_\lambda,
\]
in $U_\lambda\cap\{T<\theta\}$, where $R_\lambda$ is H\"older
continuous for $T\leq \theta$, and vanishes on $\pa\Omega$; as a
consequence, $R_\lambda=O(d^\alpha)$.

The function $\Phi_\lambda w_\lambda$ is compactly supported in
$U_\lambda$, and may be extended, by zero, to all of $\Omega$; it
is of class $C_\sharp^{2+\alpha}(\overline\Omega)$. We may
therefore consider
\[ w_1:=\sum_{\lambda\in\Lambda}\Phi_\lambda w_\lambda,\]
which is supported near $\pa\Omega$. Now, near $\pa\Omega$,
\begin{eqnarray*}
\sum_{\lambda}L(\Phi_\lambda
w_\lambda)&=&\sum_{\lambda}\Phi_\lambda L(w_\lambda)+
  2d^2\nabla\Phi_\lambda\cdot\nabla w_\lambda+
  d^2w_\lambda\Delta\Phi_\lambda
  +(4-n)w_\lambda d\nabla d\cdot\nabla\Phi_\lambda  \\
  &=& \sum_{\lambda}g\Phi_\lambda\varphi_\lambda+R'_\lambda
  = g+f,
\end{eqnarray*}
where $f=\sum_{\lambda}R'_\lambda$ has the same properties as
$R_\lambda$. It therefore suffices to solve $Lw_2=f$ when $f$ is a
H\"older continuous function which vanishes on the boundary.
\begin{lemma}\label{lem:4}
For any $f\in C^\alpha(\overline\Omega)$, there is, for $\delta$
small enough, a function $w_2\in
C_\sharp^{2+\alpha}(\overline\Omega_\delta)$ such that
\[ Lw_2=f\text{ and }w_2=O(d^\alpha)\text{ near }\pa\Omega.\]
\end{lemma}
\begin{proof}
Consider the solution $w_\ep$ of the Dirichlet problem $Lw_\ep=f$
on a domain of the form $\{\ep<d(x)<\delta\}$, with zero boundary
data. As before, $\delta$ is taken small enough to ensure that
$d\in C^{2+\alpha}(\overline\Omega_\delta)$. Schauder theory gives
$w_\ep\in C^{2+\alpha}(\{\ep\leq d(x)\leq\delta\})$. By
assumption, $|f|\leq a d^\alpha$ for some constant $a$. Let
$A>a/(\alpha+2)(n-1-\alpha)$. Since
\[
-L(d^\alpha)=d^\alpha[(\alpha+2)(n-1-\alpha)-\alpha d\Delta d],
\]
$A d(x)^\alpha$ is a super-solution if $\delta$ is small, and the
maximum principle gives us a uniform bound on $w_\ep/d^\alpha$. By
interior regularity, we obtain that, for a sequence $\ep_n\to 0$,
the $w_{\ep_n}$ converge in $C^2$, in every compact away from the
boundary, to a solution $w_2$ of $Lw_2=f$ with $w_2=O(d^\alpha)$.
Since the right-hand side $f$ is also $O(d^\alpha)$, we obtain, by
the ``type (I)'' Theorem \ref{th:FIb}, that $w_2$  of class
$C_\sharp^{1+\alpha}(\overline\Omega_\delta)$. Theorem
\ref{th:FIIa} now ensures that $w_2$ is in fact of class
$C_\sharp^{2+\alpha}(\overline\Omega_\delta)$, QED.
\end{proof}
It now suffices to take $g=-2\Delta d$ and let
\[ w_0=w_1-w_2.
\]
By construction,  $Lw_0+2\Delta d=0$ near the boundary, and $w_0$
is of class $C_\sharp^{2+\alpha}(\overline\Omega_\delta)$ if
$\delta$ is small. In addition, we know from Theorem \ref{th:L0}
that $w_1\,\rest{\pa\Omega}=(2\Delta d)/(2n-2)$, which is equal to
$-H$ on $\pa\Omega$. Lemma \ref{lem:4} gives us $w_2=O(d^\alpha)$.
We conclude that $w_0\,\rest{\pa\Omega}=-H$ on the boundary.

This completes the proof of Theorem \ref{th:III}.

% ----------------------------------------------------------------

\section{Second comparison argument and proof of Theorem \ref{th:IV}}
\label{sec:comp2}

At this stage, we have the following information, where
$\Omega_\delta=\{x : 0<d(x)<\delta\}$, for $\delta$ small enough:
\begin{enumerate}
  \item $w$ and $d\nabla w$ are bounded near $\pa\Omega$;
  \item $w=w_0+\tilde w$, where $L\tilde w=M_w(w)=O(d)$, and
  \item $w_0$
  is of class $C_\sharp^{2+\alpha}(\overline\Omega_\delta)$ for
  $\delta$ small enough.
\end{enumerate}
We wish to estimate $\tilde w$. Write $|M_w(w)|\leq cd$, where $c$
is constant.

For any constant $A>0$, define
\[ w_A:=w_0+Ad.
\]
Since $L(d)=3(2-n)d+d^2\Delta d$, we have
\[ L(w_A - w)=L(Ad-\tilde w)\leq Ad[3(2-n)+d\Delta d]+cd.
\]
Choose $\delta$ so that, say, $2(2-n)+d\Delta d\leq 0$ for $d\leq
\delta$. Then, choose $A$ so large that (i) $w_0+A\delta\geq w$
for $d=\delta$, and (ii) $(2-n)A+c\leq 0$. We then have
\[  L(w_A - w)\leq 0\text{ in }\Omega_\delta\text{ and }
      w_A - w \geq 0\text{ for }d=\delta.
\]
Next, choose $\delta$ and a constant $B$ such that
$nB+(2+Bd)\Delta d\geq 0$ on
$\Omega_\delta$. We have, by direct computation,
\[ L(d^{-2}+Bd^{-1})=-(nB+2\Delta d)d^{-1}-B\Delta d\leq 0
\]
on $\Omega_\delta$. Therefore, for any $\ep>0$,
$z_\ep:=\ep[d^{-2}+Bd^{-1}]+w_A - w$ satisfies $L z_\ep\leq 0$,
and the maximum principle ensures that $z_\ep$ has no negative
minimum in $\Omega_\delta$. Now, $z_\ep$ tends to $+\infty$ as
$d\to 0$. Therefore, $z_\ep$ is bounded below by the least value
of its negative part restricted to $d=\delta$. In other words, for
$d\leq \delta$, we have, since $w_A - w \geq 0\text{ for
}d=\delta$,
\[ w_A - w + \ep[d^{-2}+Bd^{-1}]\geq
\ep\min(\delta^{-2}+B\delta^{-1},0).
\]
Letting $\ep\to 0$, we obtain
\[ w_A - w\geq 0\text{ in }\Omega_\delta.
\]
Similarly, for suitable $\delta$ and $A$,
\[ w - w_{-A}\geq 0\text{ in }\Omega_\delta.
\]

We now know that $w$ lies between $w_0+Ad$ and $w_0-Ad$ near
$\pa\Omega$, hence $|w-w_0|=O(d)$, QED.

% ----------------------------------------------------------------

\section{Proof of main result and concluding remarks}
\label{sec:synth}

\subsection{Proof of Theorem \ref{th:V}}

At this stage, we know that
\[ L\tilde w=O(d)\text{ and }\tilde w=O(d)\]
near $\pa\Omega$ Theorem \ref{th:FIb} yields that $\tilde w$ is in
$C_\sharp^{1+\alpha}(\overline\Omega_\delta)$, for $\delta$ small
enough. It follows that $M_w(w)\in
C^{\alpha}(\overline\Omega_\delta)$. We may now use Theorem
\ref{th:FIIa} to conclude that $d^2w$ is of class $C^{2+\alpha}$
near the boundary. This completes the proof of Theorem \ref{th:V}.

Theorem \ref{th:1} now follows, as indicated in the introduction.

% ----------------------------------------------------------------

\subsection{Concluding remarks} \label{sec:concl}

The argument presented in this paper possesses some general
features, which should apply to other problems with boundary
blow-up:
\begin{itemize}
  \item Finding a singular solution $u_\Omega$ of (\ref{eq:LN})
  is equivalent to finding a bounded---in this case, classical---solution
  $w$ of an equivalent degenerate PDE.
  \item The correct regularity theorem for $w$ is stronger than
  what may be derived by scaling Schauder estimates.
  \item $w$ governs higher-order asymptotics of $u_\Omega$:
  leading-order estimates do not suffice to obtain the correct
  regularity results.
  \item The auxiliary degenerate PDE is also convenient for the
  construction of sub- and super-solutions
  which give precise control over the boundary asymptotics.
\end{itemize}
The construction of $w$ and the form of the asymptotics obey the
general rules of Fuchsian Reduction, which are not recalled (see
\cite{SK} and its references).

%%===========================================================

\end{document}